\documentclass{amsart}
\usepackage[all]{xy}
\usepackage{amsmath,amsfonts,amssymb,amsthm,epsfig,amscd,comment}

\newtheorem{Theorem}{Theorem}[section]
\newtheorem{Proposition}[Theorem]{Proposition} 
\newtheorem{Lemma}[Theorem]{Lemma}

\newtheorem{Corollary}[Theorem]{Corollary}

\theoremstyle{definition}
\newtheorem{Remark}[Theorem]{Remark}
\newtheorem{Warning}[Theorem]{Warning}

\def\A{{\mathbb{A}}}
\def\id{{\mathrm{id}}}

\def\sl{\mathfrak{sl}}

\newcommand{\Z}{\mathbb{Z}}

\newcommand{\C}{\mathbb{C}}
\newcommand{\p}{\mathbb{P}}
\newcommand{\Cone}{\operatorname{Cone}}

\def\O{{\mathcal O}}
\def\sE{{\mathcal{E}}}
\def\sL{{\mathcal{L}}}
\def\tsL{{\tilde{\mathcal{L}}}}

\def\sF{{\mathcal{F}}}
\def\sG{{\mathcal{G}}}

\def\sT{{\mathcal{T}}}
\def\sS{{\mathcal{S}}}
\def\tsT{\tilde{{\mathcal{T}}}}
\def\sP{{\mathcal{P}}}
\def\sQ{{\mathcal{Q}}}
\def\G{\mathbb{G}}
\def\wt{\widetilde}
\def\dim{\mbox{dim}}
\def\ker{\mbox{ker}}
\def\span{\mbox{span}}
\def\tB{\tilde{B}}
\def\tp{\tilde{p}}
\def\tY{\tilde{Y}}
\def\tZ{\tilde{Z}}
\def\tT{\tilde{\mathsf{T}}}
\def\ti{\tilde{i}}
\def\tj{\tilde{j}}

\newcommand{\T}{\mathsf{T}}
\newcommand{\la}{\langle}
\newcommand{\ra}{\rangle}

\newcommand{\Hom}{\operatorname{Hom}}
\newcommand{\End}{\operatorname{End}}

\newcommand{\rank}{\operatorname{rank}}

\begin{document}

\title[Equivalences and stratified flops]{Equivalences and stratified flops}

\author{Sabin Cautis}
\email{scautis@math.columbia.edu}
\address{Department of Mathematics \\ Columbia University \\ New York, NY}

\keywords{Mukai flops, stratified flops, coherent sheaves, derived equivalences, categorical actions, Cohen-Macaulay sheaves}

\begin{abstract}
We construct natural equivalences between derived categories of coherent sheaves on the local models for stratified Mukai and Atiyah flops (of type A). 
\end{abstract}
\subjclass[2000]{14E99}
\maketitle
\tableofcontents

\section{Introduction}

\subsection{Overview}

The purpose of this paper is to define in a direct way the kernels which induce equivalences between the (derived) categories of coherent sheaves on the local models for stratified Mukai or Atiyah flops (of type A). 

The relation between birational geometry and derived categories of coherent sheaves is an interesting but sometimes subtle story. On the one hand there are some general guesses of when two birational varieties should be $D$-equivalent ({\it i.e.} have isomorphic derived categories). For example, one such conjecture is that $K$-equivalence should imply $D$-equivalence. Recall that two birational varieties $Y$ and $Y'$ are $K$-equivalent if they have a common resolution $Y \xleftarrow{\pi_1} Z \xrightarrow{\pi_2} Y'$ such that $\pi_1^*(K_Y) \cong \pi_2^*(K_{Y'})$ (see \cite{H} p. 150 for details).

More generally, given a flip $X \dasharrow Y$ there should be an embedding of $DCoh(Y)$ into $DCoh(X)$. So the minimal model program can be thought of as a process which yields a birational model having the smallest derived category among all models. 

One the other hand, one can take specific birational varieties $Y$ and $Y'$ and ask for an explicit equivalence between their derived categories. One natural way to construct such an equivalence is to define a birational correspondence 
\begin{equation*}
\xymatrix{
& Z \ar[dl]_{\pi_1} \ar[dr]^{\pi_2} & \\
Y & & Y' 
}
\end{equation*}
and hope that the functor $\pi_{2*} \circ \pi_1^*: DCoh(Y) \rightarrow DCoh(Y')$ is an equivalence. Sometimes the obvious correspondence works ({\it i.e.} it induces an equivalence) and sometimes it does not (see the discussion below). 

Interesting examples of such correspondences are given by (stratified) Mukai and Atiyah flops (defined below). Such flops appear naturally in higher dimensional birational geometry. For example, Baohua Fu \cite{f} shows that two Springer maps with the same degree over a nilpotent orbit closure are connected by stratified Mukai flops. So if one can show that stratified flops induce equivalences then it follows that any two Springer resolutions are derived equivalent. 

Another example is from joint work with Joel Kamnitzer \cite{ck2} where we use explicit equivalences induced by Mukai flops to construct homological knot invariants. Stratified Mukai flops show up when one tries to generalize this construction in order to categorify the coloured $\sl_m$ Reshetikhin-Turaev knot invariants. In these cases one wants a particular equivalence -- it does not suffice to simply know the equivalence exists abstractly. 

The geometry of stratified flops has also been studied in \cite{CF} and \cite{M}. The motives and quantum cohomology rings of varieties related by flops were studied in \cite{LLW} and \cite{FW}. We will not discuss these results here although it is certainly an interesting question to relate them to the derived category of coherent sheaves. 

\subsection{Stratified Mukai and Atiyah flops}

The local model for the {\bf Mukai flop} is a correspondence which relates the cotangent bundles $T^* \p (V)$ and $T^* \p (V^\vee)$ of dual projective spaces. Notice that although $T^* \p (V)$ and $T^* \p (V^\vee)$ are isomorphic they are not naturally isomorphic (so, for example, they are not isomorphic in families). In particular, this means that the categories of coherent sheaves $Coh(T^* \p (V))$ and $Coh(T^* \p (V^\vee))$ are not naturally isomorphic. 

The correspondence is given by blowing up the zero section in $T^* \p (V)$ and then blowing down the exceptional divisor in a different direction to obtain $T^* \p (V^\vee)$. One can also blow down the zero sections of $T^* \p (V)$ and $T^* \p (V^\vee)$ to a common base $B$. To summarize, we obtain a commutative diagram 
\begin{equation*}
\xymatrix{
& W \ar[dl]_{\pi_1} \ar[dr]^{\pi_2} & \\
T^* \p (V) \ar[dr] & & T^* \p (V^\vee) \ar[dl] \\
& B. &
}
\end{equation*}

One might guess that $\pi_{2*} \circ \pi_1^*: DCoh(T^* \p (V)) \rightarrow DCoh(T^* \p (V^\vee))$ should give an equivalence. Perhaps a little surprisingly this is not the case \cite{nam1}. However, if one uses the fibre product $Z := T^* \p (V) \times_B T^* \p (V^\vee)$ instead of $W$ then the induced functor (given by pulling back to $Z$ and pushing forward) {\em does} induce a natural equivalence 
$$DCoh(T^* \p (V)) \xrightarrow{\sim} DCoh(T^* \p (V^\vee))$$
as shown by Kawamata \cite{kaw1} and Namikawa \cite{nam1}. 

Now $T^* \p (V)$ and $T^* \p (V^\vee)$ have natural one parameter deformations $\wt{T^* \p(V)}$ and $\wt{T^* \p(V^\vee)}$ over $\A^1$. The {\bf Atiyah flop} relates these two spaces by a correspondence also given by blowing up and then down the zero section. One can also blow down the zero sections to get a common space $\tilde{B}$ and obtain the analogous commutative diagram 
\begin{equation*}
\xymatrix{
& \tilde{W} \ar[dl]_{\tilde{\pi}_1} \ar[dr]^{\tilde{\pi}_2} & \\
\wt{T^* \p (V)} \ar[dr] & & \wt{T^* \p (V^\vee)} \ar[dl] \\
& \tilde{B}. &
}
\end{equation*}
This time $W$ is actually the same as the fibre product $\wt{T^* \p (V)} \times_{\tilde{B}} \wt{T^* \p (V^\vee)}$ and the functor $\tilde{\pi}_{2*} \circ \tilde{\pi}_1^*$ induces an equivalence
$$DCoh(\wt{T^* \p (V)}) \xrightarrow{\sim} DCoh(\wt{T^* \p (V^\vee)}).$$

Mukai flops and Atiyah flops have natural generalizations to cotangent bundles of Grassmannians. More precisely, the {\bf stratified Mukai flop} (of type A) relates $T^* \G(k,N)$ and $T^* \G(N-k,N)$ (the standard Mukai flop is the case $k=1$) via a diagram 
\begin{equation*}
\xymatrix{
& W \ar[dl]_{\pi_1} \ar[dr]^{\pi_2} & \\
T^* \G(k,N) \ar[dr] & & T^* \G(N-k,N) \ar[dl] \\
& B. &
}
\end{equation*}
Such flops often show up in the birational geometry of symplectic varieties. We will give a precise description of these spaces and maps in section \ref{sec1}. Perhaps surprisingly, it was shown by Namikawa in \cite{nam2} that neither $W$ nor the fibre product $Z := T^* \G(k,N) \times_B T^* \G(N-k,N)$ induce an equivalence $DCoh(T^* \G(2,4)) \xrightarrow{\sim} DCoh(T^* \G(2,4))$. At this point it is not clear what correspondence, if any, one can use to construct this equivalence. 

In \cite{ckl3} together with Joel Kamnitzer and Anthony Licata we constructed natural equivalences
$$DCoh(T^* \G(k,N)) \xrightarrow{\sim} DCoh(T^* \G(N-k,N))$$
for all $0 \le k \le N$. The construction was indirect and involved categorical $\sl_2$ actions. However, a little unexpectedly, it yielded a kernel which is a sheaf ({\it i.e.} a complex supported in one degree). In this paper we identity this sheaf in a more direct way as the pushforward of a line bundle from an open subset of the fibre product $Z$ (Theorem \ref{thm:kernel}).

There is also a deformed version of stratified Mukai flops. These are the {\bf stratified Atiyah flops} (of type A) and they relate $\wt{T^* \G(k,N)}$ and $\wt{T^* \G(N-k,N)}$ (these spaces are natural deformations of $T^* \G(k,N)$ and $T^* \G(N-k,N)$ over $\A^1$). It follows from Namikawa's work \cite{nam2} that the natural correspondence here (which is the same as the natural fibre product) does not induce an equivalence when $(k,N) = (2,4)$. In the second main theorem of this paper we use the description of the kernel used for the stratified Mukai flop (Theorem \ref{thm:kernel}) to construct another kernel (again given via pushforward of a line bundle from an open subset) which induces an equivalence 
$$DCoh(\wt{T^* \G(k,N)}) \xrightarrow{\sim} DCoh(\wt{T^* \G(N-k,N)})$$
for any $0 \le k \le N$ (Theorem \ref{thm:main}). 

In \cite{kaw2} Kawamata wrote down a tweaked version of the functor induced by $Z$ and showed that it induces an equivalence 
$$DCoh(\wt{T^* \G(2,4)}) \xrightarrow{\sim} DCoh(\wt{T^* \G(2,4)}).$$
We show that our functor agrees with his functor in this special case (Proposition \ref{prop:kawamata}).

\subsection{Outline}
In section \ref{sec1} we define all the varieties involved and their relevant deformations. We also briefly review $S_2$ sheaves and their extensions. 

In section \ref{sec2} we recall the equivalence $\T(k): DCoh(T^* \G(k,N)) \xrightarrow{\sim} DCoh(T^* \G(N-k,N))$ constructed in \cite{ckl2} and \cite{ckl3} via categorical $\sl_2$ actions. We then give a new description of the corresponding kernel $\sT(k)$ as a pushforward of a line bundle from a locally closed subset of $T^* \G(k,N) \times T^* \G(N-k,N)$ (Theorem \ref{thm:kernel}). 

In section \ref{sec3} we show that the sheaf $\sT(k)$ extends to a sheaf $\tsT(k)$ on $\wt{T^* \G(k,N)} \times_{\A^1} \wt{T^* \G(N-k,N)}$ which induces an equivalence (Theorem \ref{thm:main}). This sheaf also has a description as the pushforward of a line bundle.

In section \ref{sec4} we tie up some loose ends. First we discuss the choices of the line bundles appearing in Theorems \ref{thm:kernel} and \ref{thm:main}. Then we compute explicitly the inverse of $\tT(k)$. Finally, we explain the relation between our equivalences and the equivalence $DCoh(\wt{T^* \G(2,4)}) \xrightarrow{\sim} DCoh(\wt{T^* \G(2,4)})$ constructed by Kawamata in \cite{kaw2}.  

\subsection{Acknowledgements}

I would like to thank J\'anos Koll\'ar and S\'andor Kov\'acs for their help with questions about extensions of $S_2$ sheaves, to Rapha\"el Rouquier for some initial inspiration and to Dan Abramovich, Brendan Hassett and Joel Kamnitzer for helpful conversations. This paper was conceived while visiting MSRI -- their support is much appreciated. Research was also supported by NSF Grants 0801939/0964439.

\section{Preliminaries}\label{sec1}

In this section we discuss notation, define the main varieties we will study and review some facts regarding $S_2$ sheaves. One can skim this section on a first reading and refer back to it as a reference. 

\subsection{Notation}\label{sec:notation}

We will only consider schemes (of finite type) over $\C$. Schemes will carry some natural action of $\C^\times$ (see section \ref{sec:c*actions}). We will denote by $D(Y)$ the bounded derived category of $\C^\times$-equivariant coherent sheaves on $Y$. All maps will be $\C^\times$-equivariant. 

All functors will be derived. So, for example, if $f: Y \rightarrow X$ is a morphism then $f_*: D(Y) \rightarrow D(X)$ denotes the derived pushforward (under this convention the plain pushforward is denoted $R^0f_*$). 

\begin{Warning} There is one important exception to this rule, namely if $j: U \rightarrow Y$ is an open embedding then $j_*$ will always denote the plain (underived) pushforward. 
\end{Warning}

Recall also the terminology of Fourier-Mukai transforms (\cite{H}). Given an object $\sP \in D(X \times Y)$ (whose support is proper over $Y$) we may define the associated Fourier-Mukai transform, which is the functor
\begin{equation*}
\begin{aligned}
\Phi_\sP : D(X) &\rightarrow D(Y) \\
\sF &\mapsto {\pi_2}_* (\pi_1^* (\sF) \otimes \sP)
\end{aligned}
\end{equation*}
where $\pi_1$ and $\pi_2$ are the natural projections from $X \times Y$. The object $\sP$ is called the (Fourier-Mukai) kernel. 

Fourier-Mukai transforms have right and left adjoints which are themselves Fourier-Mukai transforms. The right and left adjoints of $ \Phi_\sP $ are the Fourier-Mukai transform with respect to the kernels 
$$ \sP_R := \sP^\vee \otimes \pi_2^* \omega_X [\dim(X)] \in D(Y \times X) \text{ and } \sP_L := \sP^\vee \otimes \pi_1^* \omega_Y [\dim(Y)] \in D(Y \times X).$$ 
If $\sP$ induces an equivalence $\Phi_{\sP}$ then its inverse is induced by $\sP_L \cong \sP_R$. 

We can express composition of Fourier-Mukai transform in terms of their kernel. If $ X, Y, Z $ are varieties and $\Phi_\sP : D(X) \rightarrow D(Y),  \Phi_\sQ : D(Y) \rightarrow D(Z) $ are Fourier-Mukai transforms, then $ \Phi_\sQ \circ \Phi_\sP $ is a Fourier-Mukai transform with respect to the kernel
\begin{equation*}
\sQ * \sP := {\pi_{13}}_*(\pi^*_{12}(\sP) \otimes \pi^*_{23}(\sQ)).
\end{equation*}
The operation $*$ is associative. Moreover by \cite{H} remark 5.11, we have $ (\sQ * \sP)_R \cong \sP_R * \sQ_R $ and $ (\sQ * \sP)_L \cong \sP_L * \sQ_L $. 

\subsection{Varieties}\label{sec:defvar}

Fix $N$ once and for all. We will write $Y(k)$ and $\tY(k)$ for $T^* \G(k,N)$ and $\wt{T^* \G(k,N)}$ respectively. If $Y$ is a variety then $D(Y)$ will denote the bounded derived category of coherent sheaves on $Y$. We usually denote closed immersions by $i$ and open immersions by $j$.  

Suppose from now on $2k \le N$. Recall the standard description of cotangent bundles to Grassmannians
$$Y(k) = \{ (X,V): X \in \End(\C^N), 0 \xrightarrow{k} V \xrightarrow{N-k} \C^N, X \C^N \subset V \text{ and } X V \subset 0 \}$$
where the arrows indicate the codimension of the inclusions. Forgetting $X$ corresponds to the projection $Y(k) \rightarrow \G(k,N)$ while forgetting $V$ gives a resolution $p(k): Y(k) \rightarrow \overline{B(k)}$ where $B(k)$ is the nilpotent orbit
\begin{equation*}
B(k) := \{ X \in \End(\C^N) : X^2 = 0 \text{ and } \dim(\ker X) = N-k \}.
\end{equation*}
The condition $\dim (\ker X) = N-k$ is equivalent to $\rank (X) = k$. 


Now we also have the resolution $p(N-k): Y(N-k) \rightarrow \overline{B(k)}$. Thus we get two resolutions
$$Y(k) \xrightarrow{p(k)} \overline{B(k)} \xleftarrow{p(N-k)} Y(N-k)$$
which is the basic example of a {\bf stratified Mukai flop} between $Y(k)$ and $Y(N-k)$. 
Consider the fibre product
\begin{eqnarray*}
Z(k) &:=& Y(k) \times_{\overline{B(k)}} Y(N-k) \\
&=& \{ 0 \overset{k}{\underset{N-k}{\rightrightarrows}} \begin{matrix} V_1 \\ V_2 \end{matrix} \overset{N-k}{\underset{k}{\rightrightarrows}} \C^N : X \C^N \subset V_1, X \C^N \subset V_2, XV_1 \subset 0, X V_2 \subset 0 \}.
\end{eqnarray*}
The scheme $Z(k)$ is equi-dimensional of dimension $2k(N-k)$. It consists of $(k+1)$ irreducible components $Z_s(k)$ ($s = 0, \dots, k)$ where
\begin{equation*}
Z_s(k) := \overline{p(k)^{-1}(B(k-s)) \times_{B(k-s)} p(N-k)^{-1}(B(k-s))}.
\end{equation*}
Notice that 
$$\{ (X,V_1,V_2) \in Z(k): \dim (\ker X) \ge N-s \} = Z_{k-s}(k) \cup Z_{k-s+1}(k) \cup \dots \cup Z_k(k)$$
$$\{ (X,V_1,V_2) \in Z(k): \dim (V_1 \cap V_2) \ge s \} = Z_{k-s}(k) \cup Z_{k-s-1}(k) \cup \dots \cup Z_0(k)$$
so we have natural increasing and decreasing filtrations of $Z(k)$. In particular, one can describe $Z_s(k)$ more directly as
$$Z_s(k) = \{(X,V_1,V_2) \in Z(k): \dim (\ker X) \ge N-k+s \text{ and } \dim (V_1 \cap V_2) \ge k-s \}.$$

Since $\span(V_1,V_2) \subset \ker X$ it follows that $\dim (\ker X) + \dim (V_1 \cap V_2) \ge N$ on $Z(k)$. We define the open subscheme
$$Z^o(k) := \{(X,V_1,V_2) \in Z(k): N+1 \ge \dim (\ker X) + \dim (V_1 \cap V_2) \} \subset Z(k)$$
and $Z^o_s(k) := Z_s(k) \cap Z^o(k)$. Abusing notation slightly, we denote all the open inclusions $j: Z^o_s(k) \rightarrow Z_s(k)$. We also denote by $i$ all the closed inclusions of $Z_s(k)$ or $Z(k)$ into $Y(k) \times Y(N-k)$. 

On $Y(k)$ we have the natural vector bundle whose fibre over $(X,V)$ is $V$. Abusing notation we denote this vector bundle by $V$. Similarly, we have vector bundles $V_1$ and $V_2$ on $Z(k)$ as well as quotient bundles $\C^N/V_1$ and $\C^N/V_2$. 

We will also use the following natural correspondences between $Y(k)$ and $Y(k+r)$
\begin{eqnarray*}
W^r(k) := \{(X,V_1,V_2) &:& X \in \End(\C^N), 0 \xrightarrow{k} V_1 \xrightarrow{r} V_2 \xrightarrow{N-k-r} \C^N, \\
& & X \C^N \subset V_1 \text{ and } XV_2 \subset 0 \}.
\end{eqnarray*}
The two projections give us an embedding $W^r(k) \subset Y(k) \times Y(k+r)$. Notice that if $r = N-2k$ then we get the standard birational correspondence
$$Y(k) \leftarrow W^{N-2k}(k) \rightarrow Y(N-k)$$
relating stratified Mukai flops. 

\subsection{Deformed varieties}

The varieties $Y(k) = T^* \G(k,N)$ have natural one parameter deformations 
$$\tY(k) = \{(X,V,x): x \in \C, X \in \End(\C^N), 0 \subset V \subset \C^N, \dim(V) = k \text{ and } \C^N \xrightarrow{X-x \cdot \id} V \xrightarrow{X+x \cdot \id} 0 \}$$
over $\A^1$. The whole picture from the previous section can be repeated for these deformed varieties. 

We can look at 
$$\tB(s) := \{X \in \End(\C^N): X^2 = x^2 \cdot \id \text{ and } \dim (\ker (X - x \cdot \id)) = N-s \}$$
and define $\tp(k): \tY(k) \rightarrow \overline{\tB(k)}$ as the map which forgets $V$ (this is a resolution). We also have the map $\tp(N-k): \tY(N-k) \rightarrow \overline{\tB(k)}$ which forgets $V$ (and maps $x$ to $-x$). As before, we get two resolutions
$$\tY(k) \xrightarrow{\tp(k)} \overline{\tB(k)} \xleftarrow{\tp(N-k)} \tY(N-k)$$
which is the basic example of a {\bf stratified Atiyah flop} between $\tY(k)$ and $\tY(N-k)$. 

The variety $Z(k)$ also deforms to 
\begin{eqnarray*}
\tZ(k) &:=& \tY(k) \times_{\overline{\tB(k)}} \tY(N-k) \\
&=& \{ 0 \overset{k}{\underset{N-k}{\rightrightarrows}} \begin{matrix} V_1 \\ V_2 \end{matrix} \overset{N-k}{\underset{k}{\rightrightarrows}} \C^N : (X-x \cdot \id) \C^N \subset V_1, (X + x \cdot \id) \C^N \subset V_2 \\
& & (X + x \cdot \id) V_1 \subset 0, (X - x \cdot \id) V_2 \subset 0 \}.
\end{eqnarray*}
However, $\tZ(k)$ is now irreducible and there are no deformed analogues of $Z_s(k)$. Notice that $\tZ(k)$ is naturally a subscheme of $\tY(k) \times_{\A^1} \tY(N-k)$ (where the second projection $\tY(N-k) \rightarrow \A^1$ maps $(X,V,x) \mapsto -x$). 

Next, we can define the open subscheme
$$\tZ^o(k) := \{(X,V_1,V_2,x) \in \tZ(k): N+1 \ge \dim (\ker(X - x \cdot \id)) + \dim (V_1 \cap V_2) \} \subset \tZ(k).$$ 
Notice that if $x \ne 0$ then $V_1$ and $V_2$ are uniquely determined by $X$ as the kernels of $X+x \cdot \id$ and $X-x \cdot \id$ respectively. In particular, this means that $\tZ^o(k)$ contains all the fibres over $x \ne 0$ since when $x \ne 0$ we have $V_1 \cap V_2 = 0$. 
We denote the open inclusion $\tj: \tZ^o(k) \rightarrow \tZ(k)$. Also, we denote by $\ti$ the closed inclusion of $\tZ(k)$ into $\tY(k) \times_{\A^1} \tY(N-k)$. As before we also have natural vector bundles $V_1,V_2,\C^N/V_1$ and $\C^N/V_2$. 

Finally, note that the varieties $W^r(k)$ do not deform. The analogue of the standard birational correspondence $Y(k) \leftarrow W^{N-2k}(k) \rightarrow Y(N-k)$ is just
$$\tY(k) \leftarrow \tZ(k) \rightarrow \tY(N-k)$$
relating stratified Atiyah flops. 

\subsection{$\C^\times$ actions}\label{sec:c*actions}

All varieties in the previous section carry a natural $\C^\times$ action. For $t \in \C^\times$, this is given by scaling $X \mapsto t^2 \cdot X \in \End(\C^N)$ and similarly (for the deformed varieties) $x \mapsto t^2 \cdot x \in \A^1$. We will always work $\C^\times$-equivariantly in this paper. 

If a variety $Y$ carries a $\C^\times$ action we denote by $\O_Y\{k\}$ the structure sheaf of $Y$ with non-trivial $\C^\times$ action of weight $k$.  More precisely, if $f \in \O_Y(U)$ is a local function then, viewed as a section $f' \in \O_Y\{k\}(U)$, we have $t \cdot f' = t^{-k}(t \cdot f)$. If $\mathcal{M}$ is a $\C^\times$-equivariant sheaf then we define $\mathcal{M}\{k\} := \mathcal{M} \otimes \O_Y \{k\}$.  Notice that this means we get a natural map of sheaves $X: \C^N \rightarrow V \{2\}$.  

\subsection{Extensions of $S_2$ sheaves}\label{sec:S2sheaves}

By a scheme we mean a separated scheme of finite type over $\C$. We review some known facts about $S_2$ sheaves. 

Let $Y$ be a scheme and $\sF$ a coherent sheaf on $Y$. Recall that if $Y$ is normal (or $S_2$) then $\sF$ is $S_2$ if and only if it is torsion free and reflexive. Equivalently, $\sF$ is $S_2$ if and only if for any open $U \subset Y$ 
$$H_X^i(U, \sF|_U) = 0$$
where $X \subset U$ is a locally closed subscheme and $0 \le i \le \mbox{min}(1, \mbox{codim}(X)-1)$. We will use both of these descriptions. 


\begin{Lemma}\label{lem:S2ext} Let $U, Y$ be schemes with $j: U \rightarrow Y$ an open embedding such that $Z := Y \setminus U \subset Y$ has codimension at least two. If $\sG$ is an $S_2$ coherent sheaf on $Y$ then $\sG \cong j_* j^* \sG$. Similarly, if $\sF$ is a coherent sheaf on $U$ then $\sF \cong j^* j_* \sF$.  
\end{Lemma}
\begin{proof}
We use the standard long exact sequence of local cohomology
$$0 \rightarrow H_Z^0(Y, \sG) \rightarrow H^0(Y, \sG) \rightarrow H^0(U, j^* \sG) \rightarrow H_Z^1(Y, \sG) \rightarrow \dots.$$
Since $\sG$ is $S_2$ and $\mbox{codim}(Z) \ge 2$ we have $H_Z^0(Y, \sG) = 0$ and $H_Z^1(Y, \sG) = 0.$ Thus $H^0(Y, \sG) \xrightarrow{\sim} H^0(U, j^* \sG).$ But we can replace $Y$ by any other open set $V \supset U$ and get the same isomorphism. Thus $\sG = j_* j^* \sG$. 

Now consider the composition of adjoint maps 
$$j_* \sF \rightarrow (j_* j^*) j_* \sF = j_* (j^* j_*) \sF \rightarrow j_* \sF.$$
The left-most map is an isomorphism by the result above while the composition is the identity by the formal properties of adjoint maps. Thus the right-most map is also an isomorphism. Now consider the exact triangle $j^* j_* \sF \xrightarrow{\phi} \sF \rightarrow \Cone(\phi)$. Applying $j_*$ we get that $j_* \phi$ is an isomorphism so $j_* \Cone(\phi) = 0$. This means $\Cone(\phi) = 0$ and hence $\phi$ is an isomorphism. 
\end{proof}

\begin{Proposition}\label{prop:S2ext} 
Let $U, Y$ be schemes with $j: U \rightarrow Y$ an open embedding such that $Y \setminus U \subset Y$ has codimension at least two. If $\sF$ is an $S_2$ coherent sheaf on $U$ then $j_* \sF$ is a coherent $S_2$ sheaf on $Y$. Moreover, it is the unique such $S_2$ sheaf whose restriction to $U$ is $\sF$. 
\end{Proposition}
\begin{proof}
Let us first assume $U$ and $Y$ are both $S_2$. 

Since $j^* \O_Y \cong \O_U$ then by Lemma \ref{lem:S2ext} we have $\O_Y = j_* j^* \O_Y \cong j_* \O_U$. Thus $j_* \sF$ is a coherent sheaf because $j_* \O_U \cong \O_Y$. 

Now a sheaf on $Y$ is $S_2$ if it is torsion-free and reflexive. Now the (underived) double dual $(j_* \sF)^{\vee \vee}$ is reflexive and torsion-free by construction. Also, since $\sF$ is $S_2$ this means that $(j_* \sF)^{\vee \vee}$ is isomorphic to $j_* \sF$ on $U$. Thus
$$(j_* \sF)^{\vee \vee} \cong j_* j^* (j_* \sF)^{\vee \vee} \cong j_* j^* j_* \sF \cong j_* \sF$$
where the first and last isomorphisms are by Lemma \ref{lem:S2ext}. 
Thus $j_* \sF$ is an $S_2$ coherent sheaf. To show uniqueness, suppose $\sG$ is another $S_2$ coherent extension of $\sF$. Then $j_* \sF \cong j_* j^* \sG \cong  \sG$ where the second isomorphism follows by Lemma \ref{lem:S2ext}. 

Finally, if $U$ is not $S_2$ then restrict to an open subscheme $U' \subset U$ whose complement has codimension at least two and repeat the argument above. If $Y$ is not $S_2$ then look at its $S_2$-ification $p: Y' \rightarrow Y$ and apply the argument above to conclude that $j'_* \sF$ is $S_2$ where $j': U \rightarrow Y'$. Then $j_* \sF = p_* j'_* \sF$ is $S_2$ since $p$ is a finite map (c.f. Proposition 5.4 of \cite{km}). 
\end{proof}

We will also need the following Lemma. 

\begin{Lemma}\label{lem:S2seq} Consider a short exact sequences of sheaves 
$$0 \rightarrow \sF_1 \rightarrow \sF \rightarrow \sF_2 \rightarrow 0$$
on an arbitrary scheme $Y$. If $\sF_1$ and $\sF_2$ are $S_2$ then $\sF$ is $S_2$. 
\end{Lemma}
\begin{proof}
Recall the description of $S_2$ from the beginning of section \ref{sec:S2sheaves}. Let $U \subset Y$ be an open subset and $X \subset U$ a locally closed subscheme. Since $\sF_1, \sF_2$ are $S_2$ this means that $H_X^i(U, \sF_j|_U) = 0$ for $j = 1,2$ and $0 \le i \le \mbox{min}(1, \mbox{codim}(X)-1)$. Using the standard long exact sequence
$$H_X^0(U, \sF_1|_U) \rightarrow H_X^0(U, \sF|_U) \rightarrow H_X^0(U, \sF_2|_U) \rightarrow H_X^1(U, \sF_1|_U) \rightarrow H_X^1(U, \sF|_U) \rightarrow H_X^1(U, \sF_2|_U) \rightarrow \dots$$
this means that $H_X^i(U, \sF|_U)  = 0$ for $0 \le i \le \mbox{min}(1, \mbox{codim}(X)-1)$ and hence $\sF$ is also $S_2$. 
\end{proof}

\section{Equivalences of stratified Mukai flops}\label{sec2}

In \cite{ckl3} we constructed a natural equivalence $\T: D(Y(k)) \xrightarrow{\sim} D(Y(N-k))$ induced by a kernel $\sT$. We begin by studying the construction of $\sT$. 

\subsection{The varieties $Z_s(k)$}

We begin with some results on $Z(k)$ and its irreducible components. It is helpful to remark that $Z(k)$ is quite singular. However, $Z^o(k)$ is quite well behaved. 

\begin{Lemma}\label{lem:Z^o} Each $Z_s(k)$ has a natural partial resolution 
$$Z'_s(k) := \{ 0 \xrightarrow{k-s} W_1 \overset{s}{\underset{N-2k+s}{\rightrightarrows}} \begin{matrix} V_1 \\ V_2 \end{matrix} \overset{N-k}{\underset{k}{\rightrightarrows}} \C^N : X \C^N \subset W_1, X V_1 \subset 0, X V_2 \subset 0 \} \xrightarrow{\pi_s} Z_s(k)$$
where $\pi_s$ forgets $W_1$. $Z_s^o(k)$ is smooth and the restriction of $\pi_s$ to the preimage of $Z_s^o(k)$ is an isomorphism. The complements of $Z_s^o(k)$ in $Z'_s(k)$ and $Z_s(k)$ have codimension at least two and four repectively. 
\end{Lemma}
\begin{proof}
$Z_s(k)$ has a natural resolution given by 
$$Z''_s(k) := \{ 0 \xrightarrow{k-s} W_1 \overset{s}{\underset{N-2k+s}{\rightrightarrows}} \begin{matrix} V_1 \\ V_2 \end{matrix} \overset{N-2k+s}{\underset{s}{\rightrightarrows}} W_2 \xrightarrow{k-s} \C^N : X \C^N \subset W_1, X W_2 \subset 0\}.$$
This variety is smooth because forgetting $V_1$ and $V_2$ gives us a $\G(s,N-2k+2s) \times \G(N-2k+s,N-2k+2s)$ fibration
$$\pi: Z''_s(k) \rightarrow \{0 \xrightarrow{k-s} W_1 \xrightarrow{N-2k+2s} W_2 \xrightarrow{k-s} \C^N: X\C^N \subset W_1, XW_2 \subset 0\}.$$
This base is actually isomorphic to the conormal bundle of the partial flag variety $\{0 \xrightarrow{k-s} W_1 \xrightarrow{N-2k+2s} W_2 \xrightarrow{k-s} \C^N \}$ embedded inside the product $\G(k-s,N) \times \G(N-k+s,N)$. 

Now $Z_s^o(k) \subset Z_s(k)$ is defined by the open condition
$$\dim (\ker X) + \dim ( V_1 \cap V_2 ) \le N+1.$$
Since $\dim (\ker X) \ge N-k+s$ and $\dim (V_1 \cap V_2) \ge k-s$ there are two possibilities: either $\dim (\ker X) = N-k+s$ or $\dim (V_1 \cap V_2) = k-s$. In the first case we can recover $W_2$ as the kernel of $X$ and $W_1$ and the image of $X$ -- this gives a local inverse of $\pi$. Similarly, in the second case we can recover $W_1$ as the intersection of $V_1$ and $V_2$ and $W_2$ as the span of $V_1$ and $V_2$. So $\pi$ is an isomorphism over $Z_s^o(k)$. In particular, since $Z''_s(k)$ is smooth this means $Z_s^o(k)$ is smooth.  

Finally, the complement $Z'_s(k) \setminus Z_s^o(k)$ is covered by three pieces:
\begin{itemize}
\item $\dim (V_1 \cap V_2) \ge k-s+2$ (where $s \ge 2$) or
\item $\dim (V_1 \cap V_2) \ge k-s+1$ and $\dim (\ker X) \ge N-k+s+1$ (where $s \ge 1$) or
\item $\dim (\ker X) \ge N-k+s+2$ (where $s \le k-1$).
\end{itemize}
The dimension of the first piece can be computed as the dimension of its resolution 
$$\{ 0 \xrightarrow{k-s} W_1 \xrightarrow{2} W_1' \overset{s-2}{\underset{N-2k+s-2}{\rightrightarrows}} \begin{matrix} V_1 \\ V_2 \end{matrix} \overset{N-2k+s}{\underset{s}{\rightrightarrows}} W_2 \xrightarrow{k-s} \C^N : X \C^N \subset W_1, X W_2 \subset 0 \}.$$
By forgetting $V_1, V_2$ and then $W_1'$ we get a map to 
\begin{equation}\label{eq:2}
\{0 \xrightarrow{k-s} W_1 \xrightarrow{N-2k+2s} W_2 \xrightarrow{k-s} \C^N: X\C^N \subset W_1, XW_2 \subset 0\}
\end{equation}
which is an iterated Grassmannian bundle with fibres 
$$\G(s-2, N-2k+2s-2), \G(N-2k+s-2,N-2k+2s-2) \text{ and } \G(2,N-2k+2s).$$
This has dimension 
$$(s-2)(N-2k+s) + (N-2k+s-2)(s) + 2(N-2k+2s-2) = (N-2k+s)(2s) - 4.$$
Now (\ref{eq:2}) is the conormal bundle of the flag variety $\mathbb{F}(k-s, N-2k+2s;N)$ inside $\G(k-s,N) \times \G(N-k+s,N)$. So its dimension is 
$$\dim(\G(k-s,N)) + \dim(\G(N-k+s,N)) = 2(N-k+s)(k-s).$$
This means that the dimension of the first piece is 
$$(N-2k+s)(2s) - 4 + 2(N-k+s)(k-s) = 2k(N-k) - 4.$$
Now the dimension of $Z_s(k)$ is $2k(N-k)$ so the codimension of the first piece is $4$.

The codimensions of the second and third pieces are computed similarly. For the second piece we use the resolution 
$$\{ 0 \xrightarrow{k-s} W_1 \xrightarrow{1} W_1' \overset{s-1}{\underset{N-2k+s-1}{\rightrightarrows}} \begin{matrix} V_1 \\ V_2 \end{matrix} \overset{N-2k+s-1}{\underset{s-1}{\rightrightarrows}} W_2 \xrightarrow{k-s+1} \C^N : X \C^N \subset W_1, X W_2 \subset 0 \}$$
which has dimension $2k(N-k)-2$ (codimension $2$). 

For the third piece we use the resolution 
$$\{ 0 \xrightarrow{k-s} W_1 \overset{s}{\underset{N-2k+s}{\rightrightarrows}} \begin{matrix} V_1 \\ V_2 \end{matrix} \overset{N-2k+s+2}{\underset{s+2}{\rightrightarrows}} W_2 \xrightarrow{k-s+2} \C^N : X \C^N \subset W_1, X W_2 \subset 0 \}$$
which has dimension $2k(N-k)-4$ again. 

To show that the codimension of $Z_s(k) \setminus Z_s^o(k) \subset Z_s(k)$ is at least four the same argument as above works except in the second case. There one needs to use the resolution 
$$\{0 \xrightarrow{k-s+1} W_1 \overset{s-1}{\underset{N-2k+s-1}{\rightrightarrows}} \begin{matrix} V_1 \\ V_2 \end{matrix} \overset{N-2k+s-1}{\underset{s-1}{\rightrightarrows}} W_2 \xrightarrow{k-s-1} \C^N : X \C^N \subset W_1, X W_2 \subset 0 \}$$
which has dimension $2k(N-k) - 2(N-2k+2s)$. This is of codimension $2(N-2k+2s) \ge 4s \ge 4$ since $N \ge 2k$ and $s \ge 1$. 
\end{proof}

\begin{Corollary}\label{cor:Z^o}
$\pi_{s*} \O_{Z'_s(k)} \cong j_* \O_{Z^o_s(k)}$.
\end{Corollary}
\begin{proof}
We can factor $j: Z^o_s(k) \rightarrow Z_s(k)$ as
$$Z^o_s(k) \xrightarrow{j'} Z'_s(k) \xrightarrow{\pi_s} Z_s(k).$$
We will see in the proof of Proposition \ref{prop:kernel} that $Z'_s(k)$ can be expressed as the intersection of the expected dimension of two smooth varieties (inside a smooth ambient variety). This means $Z'_s(k)$ is a local complete intersection and hence $S_2$ (a local complete intersection is Cohen-Macaulay which is $S_n$ where $n$ is the dimension of the variety and, in particular, $S_2$ if $n \ge 2$). 

Also, by Lemma \ref{lem:Z^o} the codimension of the complement of $Z^o_s(k)$ inside $Z'_s(k)$ is at least two. By Proposition \ref{prop:S2ext} this means that $j'_* \O_{Z^o_s(k)} \cong \O_{Z'_s(k)}$. Thus
$$\pi_{s*} \O_{Z'_s(k)} \cong \pi_{s*} j'_* \O_{Z^o_s(k)} \cong j_* \O_{Z^o_s(k)}.$$
\end{proof}

\begin{Remark}\label{Rem1}
By Lemma \ref{lem:Z^o} we can identify $Z_s^o(k)$ with 
\begin{eqnarray*}
\{ 0 \xrightarrow{k-s} W_1 \overset{s}{\underset{N-2k+s}{\rightrightarrows}} \begin{matrix} V_1 \\ V_2 \end{matrix} \overset{N-2k+s}{\underset{s}{\rightrightarrows}} W_2 \xrightarrow{k-s} \C^N &:& X \C^N \subset W_1, X W_2 \subset 0, \\
& &  \dim (\ker X) + \dim (V_1 \cap V_2) \le N+1 \}
\end{eqnarray*}
where the isomorphism to $Z_s^o(k)$ is obtained by forgetting $W_1$ and $W_2$. 
\end{Remark}

\begin{Lemma}\label{lem:intersections} Using the identification of $Z_s^o(k)$ from Remark \ref{Rem1}, the intersection
\begin{eqnarray*}
D^o_{s,+}(k) &:=& 
Z_s^o(k) \cap Z_{s+1}(k) \\
&=& \{(X,V_1,V_2) \in Z_s^o(k): \dim (\ker X) = N-k+s+1 \text{ and } \dim (V_1 \cap V_2) = k-s \}
\end{eqnarray*}
is a divisor in $Z_s^o(k)$ which is the locus where $X: \C^N/W_2 \rightarrow W_1 \{2\}$ is not an isomorphism. More precisely
\begin{equation}\label{eq:A}
\O_{Z_s^o(k)}([D^o_{s,+}(k)]) \cong \det(\C^N/W_2)^\vee \otimes \det(W_1) \{2(k-s)\}.
\end{equation}
Similarly, the intersection 
\begin{eqnarray*}
D^o_{s,-}(k) &:=& 
Z_s^o(k) \cap Z_{s-1}(k) \\
&=& \{(X,V_1,V_2) \in Z_s^o(k): \dim (\ker X) = N-k+s \text{ and } \dim (V_1 \cap V_2) = k-s+1 \}
\end{eqnarray*}
is a divisor in $Z_s^o(k)$ which is the locus where the inclusion $V_1/W_1 \rightarrow W_2/V_2$ is not an isomorphism. More precisely
\begin{equation}\label{eq:B}
\O_{Z_s^o(k)}([D^o_{s,-}(k)]) \cong \det(V_1/W_1)^\vee \otimes \det(W_2/V_2).
\end{equation}

Finally, $Z_s^o(k) \cap Z_{s'}(k) = \emptyset$ if $|s-s'| > 1$. 
\end{Lemma}
\begin{proof}
$D^o_{s,+}(k) \subset Z^o_s(k)$ is the locus where $\dim (\ker X) = N-k+s+1$. This is precisely the locus where $X: \C^N/W_2 \rightarrow W_1 \{2\}$ does not induce an isomorphism. Forgeting $X$ gives us a map from $D^o_{s,+}(k)$ to the Schubert locus $\{(V_1,V_2): \dim (V_1 \cap V_2) = k-s\}$. Since $GL(\C^N)$ acts transitively on the Schubert locus this map is a fibration. One can check that the fibres are irreducible and hence $D^o_{s,+}(k)$ is irreducible. 

Also, it is not hard to see that the scheme-theoretic locus where $X: \C^N/W_2 \rightarrow W_1 \{2\}$ does not induce an isomorphism is reduced so we conclude that 
$$\O_{Z_s^o(k)}([D^o_{s,+}(k)]) \cong \det(\C^N/W_2)^\vee \otimes \det(W_1) \{2(k-s)\}$$
This proves (\ref{eq:A}). The proof of (\ref{eq:B}) is similar. 

Finally, suppose $s'-s > 1$. The locus $Z_s(k) \cap Z_{s'}(k)$ corresponds to $\dim (\ker X) \ge N-k+s'$ and $\dim (V_1 \cap V_2) \ge k-s$. This means 
$$\dim (\ker X) + \dim (V_1 \cap V_2) \ge N+s'-s> N+1$$ 
so the intersection lies outside the open subscheme $Z^o_s(k)$ and hence $Z^o_s(k) \cap Z_{s'}(k) = \emptyset$. The case $s'-s < -1$ is dealt with similarly. 
\end{proof}

\begin{Corollary}\label{cor:intersections} On $Z_s^o(k)$ we have
$$\O_{Z_s^o(k)}([D^o_{s,+}(k)] - [D^o_{s,-}(k)]) \cong \det(\C^N/V_1)^\vee \otimes \det(V_2) \{2(k-s)\} .$$
\end{Corollary}
\begin{proof}
This is a direct consequence of Lemma \ref{lem:intersections} since 
\begin{eqnarray*}
& & \det(\C^N/W_2)^\vee \otimes \det(W_1) \otimes \det(V_1/W_1) \otimes \det(W_2/V_2)^\vee \\
&\cong& \det(\C^N/V_2)^\vee  \otimes \det(V_1) \\
&\cong& (\det(\C^N)^\vee \otimes \det(V_2)) \otimes (\det(\C^N) \otimes \det(\C^N/V_1)^\vee) \\
&\cong& \det(\C^N/V_1)^\vee \otimes \det(V_2).
\end{eqnarray*}
\end{proof}

\subsection{The kernel $\sT(k)$}

In \cite{ckl3} we defined a kernel $\sT(k) \in D(Y(k) \times Y(N-k))$ as the (right) convolution of a complex
\begin{equation}\label{eq:C}
\Theta_k(k) \rightarrow \Theta_{k-1}(k) \rightarrow \dots \rightarrow \Theta_1(k) \rightarrow \Theta_0(k).
\end{equation}
Here $\Theta_s(k) = \sF^{(N-2k+s)}(N-k) * \sE^{(s)}(k)[-s]\{s\}$ where
\begin{eqnarray*}
\sE^{(s)}(k) &=& \O_{W^{s}(k-s)} \otimes \det(\C^N/V_1)^{-s} \otimes \det(V')^{s} \{s(k-s)\} \\
\sF^{(N-2k+s)}(N-k) &=& \O_{W^{N-2k+s}(k-s)} \otimes \det(V_2/V')^{s} \{k(N-2k+s)\} 
\end{eqnarray*}
and 
\begin{eqnarray*}
W^s(k-s) = \{(X,V_1,V') &:& X \in \End(\C^N), 0 \xrightarrow{k-s} V' \xrightarrow{s} V_1 \xrightarrow{N-k} \C^N, \\
& & X \C^N \subset V' \text{ and } XV_1 \subset 0 \}.
\end{eqnarray*}
\begin{eqnarray*}
W^{N-2k+s}(k-s) = \{(X,V_1,V') &:& X \in \End(\C^N), 0 \xrightarrow{k-s} V' \xrightarrow{N-2k+s} V_2 \xrightarrow{k} \C^N, \\
& & X \C^N \subset V' \text{ and } XV_2 \subset 0 \}
\end{eqnarray*}
as in section \ref{sec:defvar} (the composition operator $*$ was defined in section \ref{sec:notation}). Notice that 
$$\sE^{(s)}(k) \in D(Y(k) \times Y(k-s)) \text{ and } \sF^{(N-2k+s)}(N-k) \in  D(Y(k-s) \times Y(N-k)).$$

\subsubsection{Convolutions} 

Let us recall the definition of a (right) convolution in a triangulated category (see \cite{GM} section IV, exercise 1).

Let $ (A_\bullet, f_\bullet) = A_n \xrightarrow{f_n} A_{n-1} \rightarrow \cdots \xrightarrow{f_1} A_0 $ be a sequence of objects and morphisms such that $ f_i \circ f_{i+1} = 0 $. Such a sequence is called a complex. A (right) convolution of a complex $ (A_\bullet, f_\bullet) $ is any object $ B $ such that there exist
\begin{enumerate}
\item objects $ A_0 = B_0, B_1, \dots, B_{n-1}, B_n = B $ and
\item morphisms $ g_i : B_i [-i] \rightarrow A_i $, $ h_i : A_i \rightarrow B_{i-1}[-(i-1)] $ (with $ h_0 = id $)
\end{enumerate}
such that
\begin{equation} \label{eq:distB}
B_i[-i] \xrightarrow{g_i} A_i \xrightarrow{h_i} B_{i-1}[-(i-1)]
\end{equation}
is a distinguished triangle for each $ i $ and $ g_{i-1} \circ h_{i} = f_i $. Such a collection of data is called a Postnikov system. Notice that in a Postnikov system we also have $f_{i+1} \circ g_i = (g_{i+1} \circ h_i) \circ g_i = 0 $ since $h_i \circ g_i = 0$.

The convolution of a complex need not exist nor is it always unique. However, in the case of the complex (\ref{eq:C}) we showed in \cite{ckl1, ckl2, ckl3} that: 

\begin{Theorem}\label{thm:ckl3}
The complex 
$$\Theta_k(k) \rightarrow \Theta_{k-1}(k) \rightarrow \dots \rightarrow \Theta_1(k) \rightarrow \Theta_0(k)$$
has a unique convolution $\sT(k) \in D(Y(k) \times Y(N-k))$ which induces an equivalence of stratified Mukai flops
$$\T(k): D(Y(k)) \xrightarrow{\sim} D(Y(N-k)).$$
\end{Theorem}

\subsubsection{Explicit description of $\sT(k)$}

Next we identify $\sT(k)$ as an $S_2$ extension of sheaves. 

\begin{Proposition}\label{prop:kernel}
We have 
$$\Theta_s(k) \cong i_* j_* \O_{Z_s^o(k)} \otimes \det(\C^N/V_1)^{-s} \otimes \det(V_2)^s [-s] \{k(N-k)-(k-s)^2 + s\}$$
as a sheaf on $Y(k) \times Y(N-k)$.
\end{Proposition}
\begin{proof}
Ignoring the $\{\cdot\}$ shift we have 
$$\Theta_s(k) \cong \pi_{13*}(\pi_{12}^* \O_{W^s(k-s)} \otimes \pi_{23}^* \O_{W^{N-2k+s}(k-s)} \otimes \det(\C^N/V_1)^{-s} \otimes \det(V')^s \otimes \det(V_2/V')^s).$$
Now 
$$\pi_{12}^{-1}(W^s(k-s)) \cap \pi_{23}^{-1}(W^{N-2k+s}(k-s)) \subset Y(k) \times Y(k-s) \times Y(N-k)$$
is the irreducible scheme 
$$Z'_s(k) = \{ 0 \xrightarrow{k-s} V' \overset{s}{\underset{N-2k+s}{\rightrightarrows}} \begin{matrix} V_1 \\ V_2 \end{matrix} \overset{N-k}{\underset{k}{\rightrightarrows}} \C^N : X \C^N \subset V', X V_1 \subset 0, X V_2 \subset 0 \}$$
which, as we saw in Lemma \ref{lem:Z^o}, is a partial resolution of $Z_s(k)$. In particular, 
$$\dim (Z'_s(k)) = 2k(N-k) = \text{ expected dimension of } \pi_{12}^{-1}(W^s(k-s)) \cap \pi_{23}^{-1}(W^{N-2k+s}(k-s)).$$ 
This means $\pi_{12}^* \O_{W^s(k-s)} \otimes \pi_{23}^* \O_{W^{N-2k+s}(k-s)} \cong \O_{Z'_s(k)}$ (one can check that the intersection is reduced by a local or cohomological calculation). 

Now $\pi_{13}$ forgets $V'$ and maps $Z'_s(k)$ to $Z_s(k)$. By Corollary \ref{cor:Z^o} $\pi_{13*} \O_{Z'_s(k)} \cong i_* j_* \O_{Z_s^o(k)}$. Thus, by the projection formula, 
$$\pi_{13*} \left( \O_{Z'_s(k)} \otimes \det(\C^N/V_1)^{-s} \otimes \det(V')^s \otimes \det(V_2/V')^s \right) \cong i_* j_* \O_{Z_s^o(k)} \otimes \det(\C^N/V_1)^{-s} \otimes \det(V_2)^s.$$

Finally, the $\{\cdot\}$ shift is 
$$s(k-s)+k(N-2k+s)+s = k(N-k)-(k-s)^2+s.$$
\end{proof}

\begin{Theorem}\label{thm:kernel}
There exists a $\C^\times$-equivariant line bundle $\sL(k)$ on $Z^o(k)$ such that $\sT(k) \cong i_* j_* \sL(k)$ where $i$ and $j$ are the natural inclusions 
$$Z^o(k) \xrightarrow{j} Z(k) \xrightarrow{i} Y(k) \times Y(N-k).$$ 
This line bundle is uniquely determined by the restrictions
$$\sL(k)|_{Z^o_s(k)} \cong \O_{Z^o_s(k)}([D^o_{s,+}(k)]) \otimes \det(\C^N/V_1)^{-s} \otimes \det(V_2)^s \{k(N-k) - (k-s)^2 + s\}.$$ 
\end{Theorem}
\begin{proof}
{\bf Step 1.} First we check that there is a unique line bundle $\sL(k)$ which restricts as above. To do this we need to check that 
$$\O_{Z^o_s(k)}([D^o_{s,+}(k)]) \otimes \det(\C^N/V_1)^{-s} \otimes \det(V_2)^s \{k(N-k) - (k-s)^2 + s\}$$
and
$$\O_{Z^o_{s-1}(k)}([D^o_{s-1,+}(k)]) \otimes \det(\C^N/V_1)^{-s+1} \otimes \det(V_2)^{s-1} \{k(N-k) - (k-s+1)^2 + s-1\}$$
agree on the overlap $D^o_{s,-} = D^o_{s-1,+}$. Cancelling line bundles and shifts this is the same as 
$$\O_{Z^o_s(k)}([D^o_{s,+}(k)]) \{2(k-s+1)\} \text{ and } \O_{Z^o_{s-1}(k)}([D^o_{s-1,+}(k)]) \otimes \det(\C^N/V_1) \otimes \det(V_2)^\vee$$
restricting to the same line bundle. Now $D^o_{s,+} \cap D^o_{s,-} = \emptyset$ so the first line bundle is just $\O_{D^o_{s,-}} \{2(k-s+1)\}$. On the other hand
$$\O_{Z^o_{s-1}(k)}([D^o_{s-1,+}(k)]) \cong \O_{Z^o_{s-1}(k)}([D^o_{s-1,-}(k)]) \otimes \det(\C^N/V_1)^\vee \otimes \det(V_2) \{2(k-s+1)\}$$
using Corollary \ref{cor:intersections}. So the second line bundle restricts to $\O_{D^o_{s,-}} \{2(k-s+1)\}$ which is the same as the first line bundle. 

{\bf Step 2.} Next we show that $\sT(k)$ is an $S_2$ sheaf. The convolution of $\Theta_*(k)$ is obtained by taking repeated cones starting on the right. If we consider the convolution of the partial complex
$$\sT_{\le s}(k) := \mbox{Conv}(\Theta_s(k) \rightarrow \dots \rightarrow \Theta_0(k))$$
then we have the exact triangle
\begin{equation}\label{eq:ses0}
\sT_{\le s}(k) \rightarrow \sT_{\le s+1} \rightarrow \Theta_{s+1}(k)[s+1].
\end{equation}
Now $\sT_{\le 0}(k) = \Theta_0(k)$ is a sheaf and by induction we can assume $\sT_{\le s}(k)$ is a sheaf. Moreover, $\Theta_{s+1}(k)[s+1]$ is an $S_2$ sheaf with reduced scheme theoretic support. So by Proposition \ref{prop:kernel} (and induction) each $\sT_{\le s+1}(k)$ is also an $S_2$ sheaf with reduced scheme theoretic support. In particular, so is $\sT(k)$. 

Since each $\sT_{\le s}(k)$ has reduced scheme theoretic support it is the pushforward of a sheaf on $Z(k)$. For the rest of the proof we will work on $Z(k)$ and to simplify notation we will think of $\sT(k)$, $\sT_{\le s}(k)$ or $\Theta_s(k)$ as lying on $Z(k)$. Now by Lemma \ref{lem:Z^o} the complement $Z(k) \setminus Z^o(k)$ has codimension at least four. So by Proposition \ref{lem:S2ext}, to show $\sT(k) \cong j_* \sL(k)$ it suffices to show $j^* \sT(k) \cong \sL(k)$. 

{\bf Step 3.} Now, by Step 1, there exists a unique line bundle $\sL_{\le s}(k)$ supported on $Z^o_{\le s}(k) := Z^o_0(k) \cup \dots \cup Z^o_s(k)$ such that its restrictions to $Z^o_{s'}(k)$ (where $0 \le s' \le s \le k$) are isomorphic to
$$\O_{Z^o_{s'}(k)}((1-\delta_{s,s'}) [D^o_{s',+}(k)]) \otimes \det(\C^N/V_1)^{-s'} \otimes \det(V_2)^{s'} \{k(N-k)-(k-s')^2+s'\}.$$
The next step is to show by induction that 
$$j^* \sT_{\le s}(k) \cong \sL_{\le s}(k). $$

The base case $s=0$ follows by Proposition \ref{prop:kernel}. To prove the induction step we look at the exact triangle             
\begin{eqnarray}\label{eq:ses}
j^* \Theta_{s+1}(k)[s] \xrightarrow{D} j^* \sT_{\le s}(k) \rightarrow j^* \sT_{\le s+1}(k) 
\end{eqnarray}
induced by (\ref{eq:ses0}). Notice that $D \ne 0$ because otherwise $\sT(k)$ would be decomposable and would subsequently not induce an equivalence. 

Now, consider the standard short exact triangle
$$\O_{Z^o_{\le s}(k)}(-[D^o_{s,+}(k)]) \xrightarrow{d} \O_{Z^o_{\le s+1}(k)} \rightarrow \O_{Z^o_{s+1}(k)}$$
where the right hand map is restriction. Tensoring by $\sL_{\le s+1}(k)$ we get 
$$\sL_{\le s+1}(k)|_{Z^o_{\le s}(k)}(-[D^o_{s,+}(k)]) \xrightarrow{d} \sL_{\le s+1}(k) \rightarrow \sL_{\le s+1}(k)|_{Z^o_{s+1}(k)}.$$
Now $\sL_{\le s+1}(k)|_{Z^o_{s+1}(k)} \cong j^* \Theta_{s+1}(k)[s+1]$ and 
$$\sL_{\le s+1}(k)|_{Z^o_{\le s}(k)} \cong \sL_{\le s}(k) \otimes \O_{Z^o_{\le s}(k)}([D^o_{s,+}(k)]) \cong j^* \sT_{\le s}(k) \otimes \O_{Z^o_{\le s}(k)}([D^o_{s,+}(k)])$$
where the second isomorphism follows by induction. Thus we end up with 
$$j^* \sT_{\le s}(k) \rightarrow \sL_{\le s+1}(k) \rightarrow j^* \Theta_{s+1}(k)[s+1]$$
and subsequently the exact triangle 
\begin{equation}\label{eq:ses2}
j^* \Theta_{s+1}(k)[s] \xrightarrow{d} j^* \sT_{\le s}(k) \rightarrow \sL_{\le s+1}(k).
\end{equation}
Notice that the map $d$ is again non-zero. Comparing (\ref{eq:ses}) and (\ref{eq:ses2}) we find that $d = D$ up to some non-zero multiple (using Lemma \ref{lem:uniquemap}). Thus $j^* \sT_{\le s+1}(k) \cong \sL_{\le s+1}$ (which completes the induction).

\begin{Lemma}\label{lem:uniquemap} We have $\Hom(j^* \Theta_{s+1}(k)[s], j^* \sT_{\le s}(k)) \cong \C$. 
\end{Lemma}
\begin{proof}
Recall that by Proposition \ref{prop:kernel} we have 
$$j^* \Theta_{s+1}(k)[s] \cong \O_{Z^o_{s+1}(k)} \otimes \det(\C^N/V_1)^{-s-1} \otimes \det(V_2)^{s+1} [-1] \{\alpha_{k,s+1}\}$$
where $\alpha_{k,s+1} = k(N-k)-(k-s-1)^2+(s+1)$. Thus we get 
\begin{eqnarray*}
&& \Hom(j^* \Theta_{s+1}(k)[s], j^* \sT_{\le s}(k))  \\
&\cong& \Hom(\O_{Z^o_{s+1}(k)} \otimes \det(\C^N/V_1)^{-s-1} \otimes \det(V_2)^{s+1}[-1] \{\alpha_{k,s+1}\}, \sL_{\le s}(k)) \\
&\cong& \Hom(\O_{Z^o_{s+1}(k)} \otimes \det(\C^N/V_1)^{-s-1} \otimes \det(V_2)^{s+1}[-1]\{\alpha_{k,s+1}\}, \sL_{\le s}(k)|_{Z^o_s(k)}) \\
&\cong& \Hom(\O_{Z^o_{s+1}(k)} \otimes \det(\C^N/V_1)^{-s-1} \otimes \det(V_2)^{s+1}[-1]\{\alpha_{k,s+1}\}, \\
& & \O_{Z^o_s(k)} \otimes \det(\C^N/V_1)^{-s} \otimes \det(V_2)^{s} \{ \alpha_{k,s} \}) \\
&\cong& \Hom(\O_{Z^o_{s+1}(k)}, \O_{Z^o_{s}(k)} \otimes \det(\C^N/V_1) \otimes \det(V_2)^\vee [1] \{-2(k-s)\}) \\
&\cong& \Hom(\O_{Z^o_{s+1}(k)},  \O_{Z^o_{s}(k)}(-[D^o_{s,+}(k)] + [D^o_{s,-}(k)] [1]) \\
&\cong& \Hom(\O_{D^o_{s,+}(k)}, \O_{D^o_{s,+}(k)}) \\
&\cong& H^0(\O_{D^o_{s,+}(k)})
\end{eqnarray*}
where the second isomorphism follows since $Z^o_{s'}(k) \cap Z^o_{s+1}(k) = \emptyset$ if $s' < s$, the third last isomorphism follows by Corollary \ref{cor:intersections} and the second last isomorphism follows by Lemma \ref{lem:union}.

Since the complement of $D^o_{s,+}(k)$ in 
\begin{eqnarray*}
D_{s,+}(k) 
&=& Z_{s}(k) \cap Z_{s+1}(k) \\
&=& \{(X,V_1,V_2) \in Z(k): \dim (\ker X) \ge N-k+s+1 \text{ and } \dim (V_1 \cap V_2) \ge k-s \}
\end{eqnarray*}
has codimension at least two it follows $H^0(\O_{D^o_{s,+}(k)}) \cong H^0(D_{s,+}(k)$. We would like to show that the only $\C^\times$-equivariant functions on $D_{s,+}(k)$ are the constants. Since $\C^\times$ acts by scaling $X \in \End(\C^N)$ the closure of every $\C^\times$ orbit intersects the locus where $X = 0$. Thus 
$$H^0(D_{s,+}(k)) \cong H^0(\{(V_1,V_2) \in \G(k,N) \times \G(N-k,N): \dim (V_1 \cap V_2) \ge k-s \}) \cong \C.$$
\end{proof}

\begin{Lemma}\label{lem:union} Consider two smooth schemes $Z_1$ and $Z_2$ intersecting along a divisor $D = Z_1 \cap Z_2$. Then $\Hom(\O_{Z_1},\O_{Z_2}(-D)[1]) \cong \Hom(\O_D,\O_D)$. 
\end{Lemma}
\begin{proof}
See, for example, Corollary 4.8 of \cite{ck1}. 
\end{proof}
\end{proof}

\section{Equivalences of stratified Atiyah flops}\label{sec3}

If we restrict to $\tZ^o(k)$, which is smooth, each $Z^o_s(k) \subset \tZ^o(k)$ defines a divisor and hence a line bundle. Notice that $Z_s(k) \subset \tZ(k)$ is in general only a divisor (it may fail to be Cartier). The following is the main result of this paper. 
Recall the natural inclusions 
$$\tZ^o(k) \xrightarrow{\tj} \tZ(k) \xrightarrow{\ti} \tY(k) \times_{\A^1} \tY(N-k).$$

\begin{Theorem}\label{thm:main}
Denote by $\tsL(k)$ the line bundle $\O_{\tZ^o(k)}(\sum_{l=0}^k \binom{l+1}{2} [Z^o_l(k)]) \{k(N-2k)\}$. Then 
$$\tsT(k) := \ti_* \tj_* \tsL(k) \in D(\tY(k) \times_{\A^1} \tY(N-k))$$ 
induces an equivalence $\tT(k): D(\tY(k)) \xrightarrow{\sim} D(\tY(N-k))$ of stratified Atiyah flops. Moreover, the restriction of $\tsT(k)$ to $Y(k) \times Y(N-k)$ is $\sT(k)$. 
\end{Theorem}
\begin{proof}
Over any non-zero $x \in \A^1$ the fibre $\tZ^o(k)_x \cong \tZ(k)_x$ maps isomorphically down to $\tY(k)_x$ and $\tY(N-k)_x$ and so $\tsT(k)_x$ induces an equivalence. So to show that $\tT(k)$ is an equivalence it suffices, by Lemma \ref{lem:technical} below, to show that $\tsT(k)$ restricts to $\sT(k)$ over the central fibre $Y(k) \times Y(N-k)$.

Consider the following commutative diagram 
\begin{equation*}
\xymatrix{
Z^o(k) \ar[r]^{j} \ar[d]^{f_1} & Z(k) \ar[r]^-{i} \ar[d]^{f_2} & Y(k) \times Y(N-k) \ar[d]^{f_3} \\
\tZ^o(k) \ar[r]^{\tj} & \tZ(k) \ar[r]^-{\ti} & \tY(k) \times_{\A^1} \tY(N-k)
}
\end{equation*}
where the $f$'s denote the natural inclusions of the central fibre. We have
\begin{eqnarray*}
f_3^* \tsT(k) &\cong& f_3^* \ti_* \tj_* \tsL(k) \cong i_* f_2^* \tj_* \tsL(k) \\
&\cong& i_* j_* j^* f_2^* \tj_* \tsL(k) \cong i_* j_* f_1^* \tj^* \tj_* \tsL(k) \cong i_* j_* f_1^* \tsL(k) 
\end{eqnarray*}
where the second isomorphism on the first row follows since $\tZ(k)$ is flat over $\A^1$ and the first isomorphism on the second row is because $\tsL(k)$ and hence $f_2^* \tj_* \tsL(k)$ is an $S_2$ sheaf (we also used Lemma \ref{lem:S2ext} which applies to both $j$ and $\tj$). Since $\sT(k) \cong i_* j_* \sL(k)$ it follows that what we really need to check is that $f_1^* \tsL(k) \cong \sL(k)$. 

As before we just need to check that the restrictions of $f_1^* \tsL(k)$ and $\sL(k)$ to each $Z_s^o(k)$ agree. This involves checking that 
$$\O_{\tZ^o(k)}(\sum_{l=0}^k \binom{l+1}{2} [Z^o_l(k)])|_{Z^o_s(k)} \cong \O_{Z^o_s(k)}([D^o_{s,+}(k)]) \otimes \det(\C^N/V_1)^{-s} \otimes \det(V_2)^s \{2ks-s^2+s\}.$$
This follows from Lemma \ref{lem:centralfibre} and the fact that if $|s-l| > 1$ then $Z^o_{l}(k) \cap Z^o_s(k) = \emptyset$ so that $\O_{\tZ^o(k)}([Z^o_{l}(k)])|_{Z^o_s(k)}$ is trivial. Putting this together we get
\begin{eqnarray*}
&& \O_{\tZ^o(k)}(\sum_{l=0}^k \binom{l+1}{2} [Z^o_l(k)])|_{Z^o_s(k)} \\
&\cong& \O_{\tZ^o(k)}(\binom{s}{2} [Z^o_{s-1}(k)] + \binom{s+1}{2} [Z^o_{s}(k)] + \binom{s+2}{2} [Z^o_{s+1}(k)])|_{Z^o_s(k)} \\
&\cong& \O_{\tZ^o(k)}(-s [Z^o_{s-1}(k)] + (s+1) [Z^o_{s+1}(k)])|_{Z^o_s(k)} \{s(s+1)\} \\
&\cong& \O_{Z^o_s(k)}(-s [D^o_{s,-}(k)] + (s+1) [D^o_{s,+}(k)]) \{s(s+1)\} \\ 
&\cong& \O_{Z_s^o(k)}([D^o_{s,+}(k)]) \otimes  \det(\C^N/V_1)^{-s} \otimes \det(V_2)^s \{2s(k-s)+s(s+1)\}
\end{eqnarray*}
where the first isomorphism uses the first fact, the second uses the second fact and the last uses Corollary \ref{cor:intersections}. 
\end{proof}

\begin{Lemma}\label{lem:technical} Let $\tY \rightarrow B$ and $\tY' \rightarrow B$ be two families of smooth varieties over some base $B$ and let $\sT \in \tY \times_B \tY'$ be a relative kernel. If $\sT|_{\tY_t \times \tY'_t}$ induces an equivalence $D(\tY_t) \xrightarrow{\sim} D(\tY'_t)$ for every $t \in B$ then $\sT$ induces an equivalence $D(\tY) \xrightarrow{\sim} D(\tY')$. 
\end{Lemma}
\begin{proof}
This follows by Proposition 3.2 of \cite{Sz}. 
\end{proof}

\begin{Lemma}\label{lem:centralfibre} 
On $\tZ^o(k)$ we have $\O_{\tZ^o(k)}(\sum_{l=0}^k [Z^o_l(k)]) \cong \O_{\tZ^o(k)} \{2\}.$ 
\end{Lemma}
\begin{proof}
The fibre of $\tZ^o(k) \rightarrow \C$ over $0 \in \C$ is $Z^o(k) = \cup_{l=0}^k Z^o_l(k)$. Since $\C^\times$ acts with weight $2$ on $\C$ it follows that $\O_{\tZ^o(k)}([Z^o(k)]) \cong \O_{\tZ^o(k)} \{2\}$ and the result follows. 
\end{proof}

\section{General Remarks}\label{sec4}

\subsection{Choice of line bundles}

The choice of line bundles $\sL(k)$ and $\tsL(k)$ in Theorems \ref{thm:kernel} and \ref{thm:main} might seem a little random. One might ask, for instance, what (if anything) goes wrong if we take $\tsL(k)$ to be the trivial line bundle. 

Namikawa showed in \cite{nam2} that in the case of $\tY(2,4)$ the trivial line bundle on $\tZ(2,4)$ (or $Z(2,4)$) fails to induce and equivalence. It is therefore likely that the trivial line bundle on $\tZ(k,N)$ also fails to induce an equivalence if $k \ne 1, N-1$. Unfortunately, we do not have a conceptual reason for this. 

On the other hand, we can tensor $\tsL(k)$ by the restriction $L$ of any line bundle from $\tZ(k)$ and then $i_* j_* (\tsL(k) \otimes L)$ will still induce an equivalence. An interesting choice is $L := \det(\C^N/V_1)^\vee \otimes \det(V_2)$ because:

\begin{Lemma}\label{lem:L} We have $L \cong \O_{\tZ^o(k)}(\sum_{l=1}^k l [Z^o_l(k)]) \{-2k\}.$
\end{Lemma}
\begin{proof}
The natural inclusion map $V_2 \rightarrow \C^N/V_1$ on $\tZ^o(k)$ is an isomorphism over any general fibre $x \ne 0$ because $V_1 \cap V_2 = 0$. So $\det(V_2) \cong \det(\C^N/V_1)$ which means that $L$ restricts to the trivial line bundle over the general fibre (here we are working non-equivariantly).

Since every line bundle on the base $\A^1$ is trivial it remains to show that 
$$L|_{Z^o(k)} \cong \O_{\tZ^o(k)}(\sum_{l=1}^k l [Z^o_l(k)]) \{-2k\}|_{Z^o(k)}.$$ 
As before, it suffices to show that 
$$L|_{Z^o_s(k)} \cong \O_{\tZ^o(k)}(\sum_{l=1}^k l [Z^o_l(k)])|_{Z^o_s(k)} \{-2k\}$$
as equivariant line bundles. Now we have 
\begin{eqnarray*}
&& \O_{\tZ^o(k)}(\sum_{l=1}^k l [Z^o_l(k)])|_{Z^o_s(k)} \{-2k\} \\
&\cong& \O_{\tZ^o(k)}((s-1)[Z^o_{s-1}(k)] + s[Z^o_s(k)] + (s+1)[Z^o_{s+1}(k)])|_{Z^o_s(k)} \{-2k\} \\
&\cong& \O_{\tZ^o(k)}(-[Z^o_{s-1}(k)] + [Z^o_{s+1}(k)])|_{Z^o_s(k)} \{2s-2k\} \\
&\cong& \O_{Z^o_s(k)}(-[D^o_{s,-}(k)] + [D^o_{s,+}(k)]) \{2(s-k)\} \\
&\cong& L|_{Z^o_s(k)} 
\end{eqnarray*}
where the first isomorphism uses that $Z^o_{s}(k) \cap Z^o_{l}(k) = \emptyset$ if $|s-l| > 1$, the second uses Lemma \ref{lem:centralfibre} while the last uses Corollary \ref{cor:intersections}.
\end{proof}

Consequently, since $\binom{l+1}{2} - l = \binom{l}{2}$, we find that:
\begin{Corollary}\label{cor:sL'} The kernel $i_* j_* \O_{\tZ^o(k)}(\sum_{l = 0}^k \binom{l}{2} [Z^o_l(k)])$ induces an equivalence of stratified Atiyah flops. 
\end{Corollary}

Notice that when $k=1$ this kernel is just $i_* j_* \O_{\tZ^o(1)} \cong \O_{\tZ(1)}$ which explains why the structure sheaf of the natural correspondence for the standard Atiyah flop {\em does} induce an equivalence (as proven in \cite{kaw1} or \cite{nam1}). 
 
\subsection{Inverses}

In this section we identify the inverse of $\T(k)$. The inverse is induced by the left (or right) adjoint of $\sT(k)$. One way to compute this is to write $\sT(k)$ as the convolution of $\Theta_* (k)$ and take left adjoints to get 
$$\sT(k)_L \cong \mbox{Conv}(\Theta_0(k)_L \rightarrow \Theta_1(k)_L \rightarrow \dots \rightarrow \Theta_k(k)_L)$$
where $\sP_L$ denotes the left adjoint of $\sP$ as defined in section \ref{sec:notation}. We have
\begin{eqnarray*}
\Theta_s(k)_L 
&=& \sE^{(s)}(k)_L * \sF^{(N-2k+s)}(N-k)_L \la s \ra \\
&\cong& \sF^{(s)}(k) \la -s(N-2k+s) \ra * \sE^{(N-2k+s)}(N-k) \la (N-2k+s)(s) \ra \la s \ra  \\
&\cong& \sF^{(s)}(k) * \sE^{(N-2k+s)}(N-k) \la s \ra 
\end{eqnarray*}
where $\la s \ra$ is short-hand for $[s]\{-s\}$. So, as before, it becomes clear that $\sT(k)_L$ is an $S_2$ sheaf which is the pushforward of some line bundle (as in Theorem \ref{thm:kernel}). To identify this line bundle we just need to work with the open subvarieties $Z^o(k)$ and $\tZ^o(k)$. 

\begin{Theorem}\label{thm:inverse} The inverse of $\tT(k)$ is induced by 
$$\tsT(k)_L \cong \ti_* \tj_* \tsL'(k) \in D(\tY(k) \times_{\A^1} \tY(N-k))$$ 
where 
$$\tsL'(k) := \O_{\tZ^o(k)}(\sum_{l=0}^k \binom{l}{2} [Z_l^o(k)]) \otimes (\det(\C^N/V_1)^{\vee} \otimes \det(V_2))^{\otimes(N-2k)} \{k(N-2k)\}.$$
The inverse of $\T(k)$ is induced by the restriction of $\tsT(k)_L$ to $Y(k) \times Y(N-k)$. 
\end{Theorem}
\begin{proof}
We will show that $\tsT(k)$ and $\ti_* \tj_* \tsL'(k)$ restricts to inverses of each other on every fibre. 

On a fibre over $x \ne 0$ this is easy to see because $\ti_* \tj_* \tsL'(k)$ restricts to the structure sheaf. This is because, as in the proof of Lemma \ref{lem:L}, $(X - x \cdot \id): V_2 \xrightarrow{\sim} \C^N/V_1$ is an isomorphism (non-equivariantly). This gives the inverse to the kernel induced by $\O_{\tZ^o(k)_x}$ because, when $x \ne 0$, $\pi_{12}^{-1}(\tZ^o(k)_x)$ and $\pi_{23}^{-1}(\tZ^o(k)_x)$ intersect transversely and the intersection maps one-to-one via $\pi_{13}$ to the diagonal in $Y(k)_x \times Y(k)_x$. Hence, if $x \ne 0$, the inverse of $\tT(k)|_x$ is $(\ti_* \tj_* \tsL'(k))|_x$. 

The case $x=0$ is more interesting. We have that 
$$\tsT(k)_L = \ti_* \tj_* (\omega_{\tZ^o(k)} \otimes \tsL(k)^\vee \otimes \pi_1^* \omega^\vee_{\tY(k)}).$$
By Proposition \ref{prop:dualizingsheaf} we have
\begin{eqnarray*}
\omega_{\tZ^o(k)} \otimes \tsL(k)^\vee \otimes \pi_1^* \omega^\vee_{\tY(k)}
&\cong& \O_{\tZ^o(k)}(\sum_{l=0}^k l^2 [Z_l^o(k)]) \otimes L^{\otimes(N-2k)} \{-2k^2-2\} \\
& & \otimes \O_{\tZ^o(k)}(\sum_{l=0}^k - \binom{l+1}{2} [Z_l^o(k)]) \{-k(N-2k)\} \otimes \pi_1^* \omega^\vee_{\tY(k)} \\
&\cong& \O_{\tZ^o(k)}(\sum_{l=0}^k \binom{l}{2} [Z_l^o(k)]) \otimes L^{\otimes(N-2k)} \otimes \pi_1^* \omega^\vee_{\tY(k)} \{-kN-2\} 
\end{eqnarray*}
where, as before, $L = \det(\C^N/V_1)^{\vee} \otimes \det(V_2)$. Now 
$$\omega_{\tY(k)}|_{Y(k)} \cong \omega_{Y(k)} \otimes \O_{Y(k)}(-[Y(k)]) \cong \O_{Y(k)} \{-2k(N-k)-2\}$$ 
where we use that $\O_{\tY(k)}([Y(k)]) \cong \O_{\tY(k)} \{2\}$ and that $\omega_{T^\star \G(k,N)} \cong \O_{T^\star \G(k,N)} \{-2k(N-k)\}$ (see Lemma 3.4 of \cite{ck3}). Hence $\tsT(k)_L$ restricts on the $x=0$ fibre to 
$$\ti_* \tj_* (\O_{\tZ^o(k)} (\sum_{l=0}^k \binom{l}{2} [Z_l^o(k)]) \otimes L^{\otimes(N-2k)} \{k(N-2k)\})|_{Y(k) \times Y(N-k)}$$
which is the same as the restriction of $\ti_* \tj_* \tsL'(k)$. This concludes the proof. 
\end{proof}

\begin{Remark} Recall that $\tsT(k)$ is induced by
$$\ti_* \tj_* \O_{\tZ^o(k)}(\sum_{l=0}^k \binom{l+1}{2} [Z^o_l(k)]) \{k(N-2k)\}$$ 
which is isomorphic to 
$$\ti_* \tj_* \O_{\tZ^o(k)}(\sum_{l=0}^k \binom{l}{2} [Z^o_l(k)]) \otimes L \{k(N-2k+2)\}.$$
On the other hand, as we just learned, $\tsT(k)_L$ is induced by 
$$\ti_* \tj_* \O_{\tZ^o(k)}(\sum_{l=0}^k \binom{l}{2} [Z^o_l(k)]) \otimes L^{N-2k} \{k(N-2k)\}.$$
Thus $\tsT(k)$ and its inverse differ by tensoring by 
$$(\det(\C^N/V_1)^\vee \otimes \det(V_2))^{\otimes(N-2k-1)} \{-2k\}$$
which is a global line bundle (and likewise $\sT(k)$ and its inverse). Even more, this line bundle is the tensor product of pullbacks of line bundles from $\tY(k)$ and $\tY(N-k)$. Thus $\sT(k)$ and its inverse differ by conjugation. This fact cannot be seen directly from the definition of $\sT(k)$ coming from categorical $\sl_2$ actions. Although this seems like a weird observation, it will help in \cite{ckl4} in constructing {\it affine} braid group actions, extending the work in \cite{ck3}. 
\end{Remark}

\begin{Proposition}\label{prop:dualizingsheaf} The dualizing sheaf of $\tZ^o(k)$ is 
$$\omega_{\tZ^o(k)} \cong \O_{\tZ^o(k)}(\sum_{l=0}^k l^2 [Z_l^o(k)]) \otimes (\det(\C^N/V_1)^{\vee} \otimes \det(V_2))^{\otimes(N-2k)} \{-2k^2-2\}.$$
\end{Proposition}
\begin{proof}
Ignoring the $\C^\times$-equivariance the right hand side restricts to the trivial line bundle on any non-central fibre $\tZ^o(k)$. This agrees with $\omega_{\tZ^o(k)_x} \cong \omega_{\tY(k)_x}$ which is trivial since $\tY(k)_x$ is obtained from $T^* \G(k,N)$ by a hyper-K\"ahler rotation. It remains to show that 
$$\omega_{Z^o(k)} \cong (\O_{\tZ^o(k)}(\sum_{l=0}^k l^2 [Z_l^o(k)] + [Z^o(k)]) \otimes L^{(N-2k)} \{-2k^2 - 2\})|_{Z^o(k)}$$
where $L = \det(\C^N/V_1)^{\vee} \otimes \det(V_2)$. Now $\O_{\tZ^o(k)}([Z^o(k)]) \cong \O_{\tZ^o(k)} \{2\}$. Also, $Z^o(k)$ is the union of $k+1$ smooth components $Z^o_s(k)$ where $Z^o_s(k) \cap Z^o_{s'}(k) = \emptyset$ if $|s-s'| > 1$. Thus 
$$\omega_{Z^o(k)}|_{Z^o_s(k)} \cong \omega_{Z^o_s(k)}([D^o_{s,-}(k)] + [D^o_{s,+}(k)])$$
so it suffices to show that 
\begin{equation}\label{eq:1}
\omega_{Z^o_s(k)}([D^o_{s,-}(k)] + [D^o_{s,+}(k)]) \cong (\O_{\tZ^o(k)}(\sum_{l=0}^k l^2 [Z_l^o(k)]) \otimes L^{\otimes(N-2k)} \{-2k^2\})|_{Z^o_s(k)}
\end{equation}
for every $s=0, \dots, k$. But $\O_{\tZ^o(k)}([Z_s^o(k)]) \cong \O_{\tZ^o(k)}(\sum_{l \ne s} - [Z_l^o(k)]) \{2\}$ so the right hand side of (\ref{eq:1}) equals 
\begin{eqnarray*}
&& \O_{Z^o_s(k)}(((s-1)^2-s^2)[D_{s,-}^o(k)] + ((s+1)^2-s^2)[D_{s,+}(k)]) \otimes L^{\otimes(N-2k)} \{-2k^2 + 2s^2\} \\
&\cong& \O_{Z^o_s(k)}((-2s+1)[D^o_{s,-}(k)] + (2s+1)[D^o_{s,+}(k)]) \otimes L^{\otimes(N-2k)} \{-2k^2 + 2s^2\} \\
&\cong& \O_{Z^o_s(k)}([D^o_{s,-}(k)] + [D^o_{s,+}(k)]) \otimes L^{\otimes(N-2k+2s)} \{4s(k-s) - 2k^2 + 2s^2\}
\end{eqnarray*}
where we used Corollary \ref{cor:intersections} to get the last isomorphism. By Lemma \ref{lem:dualizingsheaf} this equals the left hand side of (\ref{eq:1}). 
\end{proof}

\begin{Lemma}\label{lem:dualizingsheaf} The dualizing sheaf of $Z^o_s(k)$ is 
$$\omega_{Z^o_s(k)} \cong (\det(\C^N/V_1)^{\vee} \otimes \det(V_2))^{\otimes(N-2k+2s)} \{-2(k-s)^2\}.$$
\end{Lemma}
\begin{proof}
As we saw in the proof of Lemma \ref{lem:Z^o} $Z_s(k)$ has a natural resolution 
$$Z''_s(k) := \{ 0 \xrightarrow{k-s} W_1 \overset{s}{\underset{N-2k+s}{\rightrightarrows}} \begin{matrix} V_1 \\ V_2 \end{matrix} \overset{N-2k+s}{\underset{s}{\rightrightarrows}} W_2 \xrightarrow{k-s} \C^N : X \C^N \subset W_1, X W_2 \subset 0\}$$
where the map $\pi: Z''_s(k) \rightarrow Z_s(k)$ is an isomorphism over $Z_s^o(k)$. Hence $\omega_{Z^o_s(k)} \cong \omega_{Z''_s(k)}|_{Z^o_s(k)}$. 

Now consider the projection map 
$$p: Z''_s(k) \rightarrow \{0 \xrightarrow{k-s} W_1 \xrightarrow{N-2k+2s} W_2 \xrightarrow{k-s} \C^N:  X \C^N \subset W_1, X W_2 \subset 0\}$$
given by forgetting $V_1$ and $V_2$. This is a $\G(s, W_2/W_1) \times \G(N-2k+s,W_2/W_1)$ fibration so the relative cotangent bundle is
\begin{eqnarray*}
\omega_p 
&\cong& (\det(V_1/W_1)^{N-2k+s} \otimes \det(W_2/V_1)^{-s}) \otimes (\det(V_2/W_1)^{s} \otimes \det(W_2/V_2)^{-N+2k-s}) \\
&\cong& (\det(\C^N/V_1)^\vee \otimes \det(V_2))^{N-2k+2s} \otimes (\det(\C^N/W_2)^\vee \otimes \det(W_1))^{-N+2k-2s}.
\end{eqnarray*}
On the other hand, 
$$W^{N-2k+2s}(k-s) \subset \{0 \xrightarrow{k-s} W_1 \xrightarrow{N-2k+2s} W_2 \xrightarrow{k-s} \C^N:  X W_2 \subset 0\}$$
is carved out by the section $X: \C^N/W_2 \rightarrow W_2/W_1 \{2\}$. Then
\begin{eqnarray*}
\omega_{W^{N-2k+2s}(k-s)} 
&\cong& \omega_{T^\star \G(N-k+s,N)} \otimes (\det(W_1)^{N-2k+2s} \otimes \det(W_2/W_1)^{-k+s}) \\
& & \otimes (\det(\C^N/W_2)^{-N+2k-2s} \otimes \det(W_2/W_1\{2\})^{k-s}) \\
&\cong& (\det(\C^N/W_2)^\vee \otimes \det(W_1))^{N-2k+2s} \\
& & \{-2(N-k+s)(k-s) + 2(k-s)(N-2k+2s) \}
\end{eqnarray*}
where we use that $\omega_{T^\star \G(a,N)} \cong \O_{T^\star \G(a,N)} \{-2a(N-a)\}$ (see Lemma 3.4 of \cite{ck3}). Thus 
$$\omega_{Z''_s(k)} \cong \omega_p \otimes \omega_{W^{N-2k+2s}(k-s)} \cong (\det(\C^N/V_1)^\vee \otimes \det(V_2))^{N-2k+2s} \{-2(k-s)^2\}$$ 
and the result follows. 
\end{proof}

\subsection{Kawamata's equivalence for $T^* \G(2,4)$}

For the rest of the paper we will ignore the $\C^\times$-equivariance. In \cite{kaw2} Kawamata studies stratified Atiyah flops for $\tY(2) = \wt{T^* \G(2,N)}$. Incidentally, Kawamata also calls these stratified Mukai flops (the terminology of stratified Atiyah flops is from \cite{nam2}). 

Kawamata constructs two functors $\Psi, \Phi: D(\tY(2)) \rightarrow D(\tY(N-2))$ which are adjoints of each other. He shows (\cite{kaw2} Theorem 3.5) that when $N=4$ these functors induce equivalences. We would like to briefly outline his construction and explain why:

\begin{Proposition}\label{prop:kawamata} The kernel $\sS(2) = i_* j_* \O_{\tZ^o(2)}(\sum_{l = 0}^2 \binom{l}{2} [Z^o_l(2)])$ and its (left or right) adjoint induce Kawamata's functors $\Psi, \Phi: D(\tY(2)) \rightarrow D(\tY(N-2))$. 
\end{Proposition}

\begin{Remark}
By Corollary \ref{cor:sL'} the kernel $\sS(2)$ induces an equivalence. Hence Kawamata's functors $\Psi$ and $\Phi$ are equivalences for all $N$. 
\end{Remark}

Inside 
$$\tY(2) = \{(X,V,x): x \in \C, X \in \End(\C^N), 0 \subset V \subset \C^N, \dim(V) = 2 \text{ and } \C^N \xrightarrow{X-x \cdot \id} V \xrightarrow{X+x \cdot \id} 0 \}$$
there are two natural strata: namely the locus where $X=0$ (isomorphic to $\G(2,N)$) and the locus where $\rank X \le 1$ (the locus $\rank X \le 2$ is the whole central fibre $Y(2)$). Kawamata first blows up the first locus and then the strict transform of the second locus to obtain
$$\tY''(2) \xrightarrow{f_1} \tY'(2) \xrightarrow{f_2} \tY(2).$$
Inside $\tY''(2)$ we denote by $E_1$ the exceptional divisor of $f_1$ and by $E_2$ the strict transform of the exceptional divisor of $f_2$. Warning: our labeling of divisors does not match precisely that in \cite{kaw2} for reasons that will soon become clear.

Kawamata then blows up $\tY(N-2)$ in the same way to obtain $\tY''(N-2)$ and identifies this smooth variety with $\tY''(2)$. To summarize, we arrive at the following commutative diagram
\begin{equation}\label{eq:kawamata}
\xymatrix{
& \tY''(2) \cong \tY''(N-2) \ar[dl]_{f} \ar[dr]^{f^+} \ar[d]^{\pi} & \\
\tY(2) \ar[rd] & \tZ(2) \ar[l]_{\pi_1} \ar[r]^{\pi_2} & \tY(N-2) \ar[ld] \\
& \overline{\tB(2)}. &
}
\end{equation}
Notice that $\tZ(2) = \tY(2) \times_{\overline{\tB(2)}} \tY(N-2)$ so the map $\pi$ exists by the universal property of fibre products. The functors $\Psi$ and $\Phi$ are then defined by 
\begin{eqnarray*}
\Psi(\cdot) &:=& f_*(f^{+*}(\cdot) \otimes \O_{\tY''(2)}([E_2])) \\
\Phi(\cdot) &:=& f^+_*(f^{*}(\cdot) \otimes \O_{\tY''(2)}((2N-5)[E_2] + (N-3)[E_1])).
\end{eqnarray*}

Notice that $Z(2)$ contains three components, namely $Z_s(2)$ where $s=0,1,2$. Similarly, $Y''(2)$ (the preimage of $Y(2)$ under $f$) also contains three components, namely $E_1, E_2$ and the proper transform of $Y(2)$ inside $\tY''(2)$ (which we denote $E_0$). The map $\pi$ then maps $E_s$ to $Z_s(2)$ for $s=0,1,2$ (the map is generically one-to-one). 

From this we see that the complement of $\tY''(2)^o := \pi^{-1}(\tZ^o(2))$ in $\tY''(2)$ has codimension at least two (otherwise it would contain a divisor whose restriction to $Y''(2)$ would be another component). Now denote by $j': \tY''(2)^o \rightarrow \tY''(2)$ the natural open immersion. Thus for any $a,b \in \Z$ 
$$j'_*(\O_{\tY''(2)^o}(a [E_1] + b [E_2])) \cong \O_{\tY''(2)}(a [E_1] + b [E_2]).$$

Using the commutative diagram 
\begin{equation*}
\xymatrix{
\tY''(2)^o \ar[d]^{\pi_o} \ar[r]^{j'} & \tY''(2) \ar[d]^{\pi} \\
\tZ^o(2) \ar[r]^{j} & \tZ(2)
}
\end{equation*} 
we find that 
\begin{eqnarray*}
R^0 \pi_* \O_{\tY''(2)}([E_2]) &\cong& R^0 \pi_* \circ j'_*(\O_{\tY''(2)^o}([E_2])) \cong j_* \circ R^0 \pi_{o*} (\O_{\tY''(2)^o}([E_2])) \\
&\cong& j_* \O_{\tZ^o(2)}([Z_2^o(2)]) \cong \sS(2)
\end{eqnarray*}
where the third isomorphism follows by the projection isomorphism since $E_2 = \pi_o^{-1} Z_2^o(2)$ and $\pi_o$ is a proper map between smooth varieties. Now 
\begin{eqnarray*}
\Psi(\cdot) 
&\cong& f_*(f^{+*}(\cdot) \otimes \O_{\tY''(2)}([E_2])) \\
&\cong& (\pi_{1*} \pi_*)(\pi^* \pi_2^*(\cdot) \otimes \O_{\tY''(2)}([E_2])) \\
&\cong& \pi_{1*} (\pi_2^*(\cdot) \otimes \pi_*(\O_{\tY''(2)}([E_2]))) \\
&\cong& \pi_{1*} (\pi_2^*(\cdot) \otimes \sS(2))
\end{eqnarray*}
which proves Proposition \ref{prop:kawamata}. 

As a consistancy check one can similarly show that 
\begin{eqnarray*}
& & R^0 \pi_* \O_{\tY''(2)}((2N-5)[E_2] + (N-3)[E_1]) \\
&\cong& j_* \O_{\tZ^o(2)}((2N-5)[Z_2^o(2)] + (N-3)[Z_1^o(2)]) \\
&\cong& j_* \O_{\tZ^o(2)}(\sum_{l=0}^2 \binom{l}{2} [Z_s^o(2)]) \otimes (\det(\C^N/V_1)^\vee \otimes \det(V_2))^{N-3} \\
&\cong& \sS(2)_L
\end{eqnarray*}
where we used that $\det(\C^N/V_1)^\vee \otimes \det(V_2) \cong \O_{\tZ^o(2)}(\sum_{l=0}^2 l [Z_l^o(2)])$. This means that $\sS(2)_L$ induces $\Phi$. Note that we already knew this because $\Psi$ and $\Phi$ were shown to be adjoints in \cite{kaw2}. 

\begin{Remark}
This same argument can be generalized from $k=2$ to arbitrary $k$. More precisely, one can iteratively blow up $\tY(k)$ starting with the smallest stratum to obtain 
$$\tY^{(k)}(k) \xrightarrow{f_1} \tY^{(k-1)}(k) \xrightarrow{f_2} \dots \xrightarrow{f_k} \tY(k).$$
We denote by $E_s$ the (proper transform of the) exceptional divisor of $f_s$. 

Then there is an isomorphism $\tY^{(k)}(k) \cong \tY^{(k)}(N-k)$ (see \cite{M}) and a commutative diagram like in (\ref{eq:kawamata}). The map $\pi: \tY^{(k)}(k) \rightarrow \tZ(k)$ identifies the divisors $E_s$ with $Z_s(k)$ for $s=0, \dots, k$. Subsequently:

\begin{Proposition}\label{prop:kawamata2} 
$$R^0 \pi_*(\O_{\tY^{(k)}(k)}(\sum_{l = 0}^k \binom{l}{2} [E_l])) \cong j_* \O_{\tZ^o(k)}(\sum_{l = 0}^k \binom{l}{2} [Z^o_l(k)])$$ 
and thus both kernels induce equivalences $D(\tY(k)) \xrightarrow{\sim} D(\tY(N-k))$. 
\end{Proposition}
\end{Remark}

\end{document}